\documentclass[12pt]{amsart}
\usepackage{amssymb}

\newcommand{\tPa}[6]{\bibitem{#1} {#2}, \emph{#3}, {#4}, to appear (#5 pages).}
\newcommand{\sPa}[5]{\bibitem{#1} {#2}, \emph{#3}, {#4}, submitted.}
\newcommand{\FinSeqs}[1]{{{#1}^{<\alephes}}}
\newcommand{\rest}[1]{\mid_{#1}}
\newcommand{\Bgp}{{\Z^\N}}

\long\def\forget#1\forgotten{}
\newcommand{\issuenumber}{21}
\newcommand{\issuemonth}{June}
\newcommand{\issueyear}{2007}

\newtheorem{thm}{Theorem}[section]
\newtheorem{prob}[thm]{Problem}

\newtheorem{issue}{Issue}

\theoremstyle{definition}

\theoremstyle{remark}

\newcommand{\alephes}{{\aleph_0}}
\newcommand{\ed}{
\newpage

\section{Unsolved problems from earlier issues}

\begin{issue}
Is $\binom{\Omega}{\Gamma}=\binom{\Omega}{\Tau}$?
\end{issue}

\begin{issue}
Is $\ufin(\cO,\Omega)=\sfin(\Gamma,\Omega)$?
And if not, does $\ufin(\cO,\Gamma)$ imply
$\sfin(\Gamma,\Omega)$?
\end{issue}

\stepcounter{issue}

\begin{issue}
Does $\sone(\Omega,\Tau)$ imply $\ufin(\Gamma,\Gamma)$?
\end{issue}

\begin{issue}
Is $\fp=\fp^*$? (See the definition of $\fp^*$ in that issue.)
\end{issue}

\begin{issue}
Does there exist (in ZFC) an uncountable set satisfying $\sone(\BG,\B)$?
\end{issue}

\stepcounter{issue}

\begin{issue}
Does $X \nin \NON(\M)$ and $Y\nin\mathsf{D}$ imply that
$X\cup Y\nin \COF(\M)$?
\end{issue}

\begin{issue}[CH]
Is $\split(\Lambda,\Lambda)$ preserved under finite unions?
\end{issue}

\begin{issue}
Is $\cov(\M)=\fo$? (See the definition of $\fo$ in that issue.)
\end{issue}

\begin{issue}
Does $\sone(\Gamma,\Gamma)$ always contain an element of cardinality $\fb$?
\end{issue}

\begin{issue}
Could there be a Baire metric space $M$ of weight $\aleph_1$ and a partition
$\mathcal{U}$ of $M$ into $\aleph_1$ meager sets where for each ${\mathcal U}'\subset\mathcal U$,
$\bigcup {\mathcal U}'$ has the Baire property in $M$?
\end{issue}

\stepcounter{issue} 

\begin{issue}
Does there exist (in ZFC) a set of reals $X$ of cardinality $\fd$ such that all
finite powers of $X$ have Menger's property $\sfin(\cO,\cO)$?
\end{issue}

\begin{issue}
Can a Borel non-$\sigma$-compact group be generated by a Hurewicz subspace?
\end{issue}

\begin{issue}[MA]
Is there an uncountable $X\sbst\R$ satisfying $\sone(\BO,\BG)$?
\end{issue}

\begin{issue}[CH]
Is there a totally imperfect $X$ satisfying $\ufin(\cO,\Gamma)$
that can be mapped continuously onto $\Cantor$?
\end{issue}

\begin{issue}[CH]
Is there a Hurewicz $X$ such that $X^2$ is Menger but not Hurewicz?
\end{issue}

\begin{issue}
Does the Pytkeev property of $C_p(X)$ imply the Menger property of $X$?
\end{issue}

\begin{issue}
Does every hereditarily Hurewicz space satisfy $\sone(\BG,\BG)$?
\end{issue}

\begin{issue}[CH]
Is there a Rothberger-bounded $G\le\Bgp$ such that $G^2$ is not Menger-bounded?
\end{issue}

\general\end{document}}

\newcommand{\Cantor}{{\{0,1\}^\N}}
\newcommand{\oo}{\infty}

\newcommand{\fb}{\mathfrak{b}}

\newcommand{\fd}{\mathfrak{d}}
\newcommand{\fp}{\mathfrak{p}}

\newcommand{\NON}{{\mathsf   {NON}}}
\newcommand{\COF}{{\mathsf   {COF}}}

\newcommand{\M}{\mathcal{M}}

\newcommand{\cov}{\mathsf{cov}}

\newcommand{\R}{\mathbb{R}}

\newcommand{\fo}{\mathfrak{od}}

\renewcommand{\split}{\mathsf{Split}}
\newcommand{\bq}{\begin{quote}}
\newcommand{\eq}{\end{quote}}
\newcommand{\cO}{\mathcal{O}}
\newcommand{\B}{\mathcal{B}}
\newcommand{\BG}{\B_\Gamma}

\newcommand{\BO}{\B_\Omega}

\newcommand{\sone}{\mathsf{S}_1}    \newcommand{\sfin}{\mathsf{S}_{fin}}

\newcommand{\ufin}{\mathsf{U}_{fin}}

\newcommand{\nin}{\not\in}

\newcommand{\NN}{{\N^\N}}
\newcommand{\N}{\mathbb{N}}
\newcommand{\Z}{\mathbb{Z}}
\newcommand{\sm}{\setminus}
\newcommand{\sbst}{\subseteq}
\newcommand{\by}[2]{\par\hfill\emph{#1}, #2}
\newcommand{\nby}[1]{\par\hfill\emph{#1}}
\newcommand{\Tau}{\mathrm{T}}
\newcommand{\CE}{\textsc{CE}}

\newcommand{\be}{\begin{enumerate}}
\newcommand{\ee}{\end{enumerate}}
\newcommand{\bi}{\begin{itemize}}
\newcommand{\ei}{\end{itemize}}
\newcommand{\itm}{\item}

\newcommand{\general}{\small\vfill\par\noindent\hrulefill\par
\noindent\textbf{Previous issues.} The previous issues of this
bulletin, which contain general information (first issue), basic
definitions, research announcements, and open problems (all
issues) are available online,
at \texttt{http://front.math.ucdavis.edu/search?\&t=\%22SPM+Bulletin\%22}
\\[0.1cm]
\textbf{Contributions.}
Please submit your contributions (announcements, discussions, and open problems)
by e-mailing us. It is preferred to write them
in \LaTeX{}.
The authors are urged to use as standard notation as possible, or otherwise give
the definitions or a reference to where the notation is explained.
Contributions to this bulletin would not require any transfer of copyright,
and material presented here can be published elsewhere.\\[0.1cm]
\textbf{Subscription.}
To receive this bulletin (free) to your
e-mailbox, e-mail us:\\
{boaz.tsaban@weizmann.ac.il}
}

\newcommand{\nArxPaper}[5]{\subsection{#2}{#4}\par\hfill{\arx{#1}}\par\hfill\emph{#3}}

\newcommand{\nAMSPaper}[4]{\subsection{#2}{#4}\par\hfill{\texttt{#1}}\par\hfill\emph{#3}}

\newcommand{\arx}[1]{\texttt{http://arxiv.org/abs/#1}}
\newcommand{\url}[1]{\bq\texttt{#1}\eq}
\newcommand{\online}[1]{The paper is available online at \url{#1}}

\title[$\mathcal{SPM}$ Bulletin \textbf{\issuenumber} (\issuemonth{} \issueyear)]{%
$\mathcal{SPM}$ Bulletin\\[0.5cm]
Issue number \issuenumber: \issuemonth{} \issueyear{} \CE{}}

\begin{document}
\maketitle

\tableofcontents

\section{Editor's note}

This is the first issue after the successful SPM07 meeting.
Lively personal impressions of two major participants appear
in Section \ref{impr}: A detailed account by Zdomskyy, including
a brief description of the talks, and a concise one by Di Maio.
I hope that this will be some sort of compensation for readers who
unfortunately could not make it to the meeting, and give us a
good reason to look forward to the next meetings.

The Workshop's official web-page is
\url{http://www.pmf.ni.ac.yu/spm2007/index.html}
The presentations are available at
\url{http://www.cs.biu.ac.il/~tsaban/SPMC07/slides.html}
Some pictures are available at
\url{http://www.thesamet.com/spm/photos/main.php}
(more pictures will be added there later).

\medskip

\by{Boaz Tsaban}{boaz.tsaban@weizmann.ac.il}

\hfill \texttt{http://www.cs.biu.ac.il/\~{}tsaban}

\section{Personal impressions from the SPM07 meeting}\label{impr}

The atmosphere of the workshop was friendly and stimulating. The
rather warm April for central Europe, and the beautiful nature
surrounding us added to the positive atmosphere. A welcome party
and a banquet were organized in a typical Slavic manner. The
participants could enjoy visiting centers of Serbian culture. This
was the background of the \textbf{III Workshop on Coverings,
Selections and Games in Topology}.

The conference was opened by the talk of \emph{V. Fedorchuk}, who
presented some results from dimension theory. In particular,
\emph{wid = w-$\infty$-C}, and hence all the properties
\emph{w-$m$-C} coincide, where $m\in\N$. Earth is really compact:
Fedorchuk is my mathematical great grandfather.

By its spirit, Fedorchuk's talk was close to that of
\emph{L.~Babinkostova}. These talks demonstrated that selection
principles is not only the ``zero-dimensional'' part of topology,
but also has natural connections with dimension theory, and there
are promising directions here.

\emph{M. Scheepers} showed in his talk how many faces the
Sakai property $S_1(\Omega,\Omega)$ has. The talk was very
impressive: one of the participants told me that he would be happy
to have this talk at least one hour longer.

The next speaker, \emph{G. Di Maio} presented a very intensive
talk. It was shown how natural it is to bring the notion of
abstract boundedness into selection principles. He also showed
connections with Arkhangel'ski-Ko\v{c}inac $\alpha_i$-properties,
and Ramsey theory. And he a musician as well!

\emph{M. Mr\v{s}evi\'{c}} presented some current results
concerning selection properties of topological (hyper)spaces
defined by means of \v{C}ech closure operators. Except for the
mathematical content of her talk, the participants could enjoy
(and join) her greetings to Ljubi\v{s}a Ko\v{c}inac on the
occasion of his 60th (he looks younger) anniversary.

Some heavy and very interesting stuff about  universal elements in
some classes of mappings and classes of $G$-spaces was presented
by \emph{S.~Iliadis}.  The number of definitions and theorems per
minute was so big, that one should be in good shape to follow them
all. Nonetheless, the main results were presented in a clear way.

\emph{P.~Blagojevi\v{c}, T.-L.~Costache, and  M.~Joita's} talks
returned me to the reality that the mathematics is very rich and
diverse. The main ideas were presented in a very gentle way, which
made it very pleasant to listen  to these talks. Blagojevi\v{c}
started with a nice riddle, and proceeded to obtain a powerful
Borsuk-Ulam type theorem. Costache and Joita presented a number of
results about $C^*$-algebras.

\emph{L.~Ko\v{c}inac} has also demonstrated the quick growth of
field of selection principles: Prior to his talk, I heard the
terms ``regularly varying, slowly varying'' sequences only from
people doing complex analysis. Some of the presented results are
nice: one can  have a feeling that it is possible to prove them
``on fingers'', but gives up after delicate places are
encountered.

In his talk, \emph{D.~Georgiou} mainly considered the different
topologies on the spaces of continuous functions and compared
there properties. Besides the huge variety of results, he did an
important job: you can find a quite long list of many classical
and important works determining  directions in topology in his
presentation.

\emph{B.~Tsaban and N.~Samet} presented a nice one-hour talk
consisting of two parts devoted to Ramsey theory. The first part
(by Tsaban) contained motivation, all needed definitions and basic
results in Ramsey theory. One could enjoy nice presentation
accompanied by clear explanations, which resulted in a  smooth and
accessible introduction to the topic.

The heavy staff was presented by Samet. He made a (successful, to
my opinion) effort to bring the participants to the front line of
Ramsey Theory and presented some delicate proofs.

The most impressive effort to force the participants to understand
a tricky and difficult proof was made by \emph{M.~Machura}. He
spoke about his result (joint with S. Shelah and B. Tsaban)
concerning finite powers of $o$-bounded groups. The importance of
the results justifies the difficulty of the proof. Macura's
didactic presentation was a work of art.

\emph{M.~Sakai} mainly spoke about local properties of spaces of
continuous functions. Some of the results were quite surprising:
It was brave of him to dare (and succeed) proving that $C_p(X)$
has the property $(\sharp)$ for every $X$.

Some results concerning selection principles in relator spaces
were presented by \emph{D.~Kocev}. He demonstrated that many
well-known facts in the realm of topological spaces can be
extended to realtor spaces. Therefore, these results do not even
use that the underlying object is a topological space in a full
strength, let alone some separation axioms, etc.

\emph{V. Vuksanovi\v{c}} gave an overview of the present status of
so-called canonical Ramsey theory. In particular, he spoke about
generalizations of such notions like the Ellentuck topology, and
mentioned some prominent results  of Todorcevic and others whose
proofs use these generalizations.

Various topologies on the set $C(X)$ of all continuous real-valued
functions were considered by \emph{V. Pavlovi\'c}. Besides the
compact-open one, he considered  topologies obtained when a
function is identified with its graph, which is treated as a point
of  some hyperspace. Analogies with the classical results were
presented. Typical to his good sense of humor, the end of the talk
was really funny.

The following  striking characterization of $I$-favorable spaces
was presented by \emph{S.~Plewik}: A compact space $X$ is
$I$-favorable if, and only if, $X$ can be represented as a limit
of $\sigma$-complete inverse system of compact metrizable spaces
with skeletal bonding maps.  By its spirit, this result is similar
to Shchepin's theory of openly generated spaces.

The talk of \emph{P.~Kalemba} was devoted to doughnuts, which are
natural generalizations of basic open sets in the Ellentuck
topology. In particular, some modifications of the \emph{base
matrix tree lemma}  were presented, and they were used   to
establish a number of inequalities between cardinal
characteristics of doughnuts.

\emph{L.~Bukovsk\'y} spoke about properties of spaces of
continuous functions defined with help of the quasi-normal
convergence: $\mathit{wQN}$-, $\mathit{QN}$-spaces, and their
modifications. One could enjoy the clear presentation of his
results by diagrams.

I talked about joint results with B. Tsaban, characterizing
several classes of hereditarily Hurewicz spaces. In the workshop
we learned that the one of the main problems from Bukovsky's talk is
solved by our characterizations. This is yet another demonstration of the
importance of such meetings.

The future of the $\mathcal{SPM}$-Bulletin was discussed during
the SPM Forum. One of the most interesting (to my opinion)
suggestions was made by S.~Plewik, who proposed  to establish a
collaboration with the AMS in order to improve the
MathSciNet-reviews of the SPM-papers.

Finally, participants had a nice opportunity to learn many
interesting facts about the life and work of  \emph{D.~Kurepa}, in
a Round Table dedicated to his centennial. Besides the information
available online, I was impressed by the personal recollections of
L.~Bukovsk\'y concerning his meetings with Kurepa.

\nby{Lyubomyr Zdomskyy}

\bigskip

\noindent
The SPM 07 workshop was a great meeting from several points of view:
\be
\itm Carefully organized (and I believe that this is not easy, especially in Serbia);
\itm The atmosphere was friendly;
\itm The level of talks was high;
\itm I feel at home in the SPM group.
\ee

\nby{Giuseppe Di Maio}

\section{Research announcements}

\nArxPaper{0704.1884}
{Coloring ordinals by reals}
{J\"org Brendle and Saka\'e Fuchino}
{We study combinatorial principles we call Homogeneity Principle $HP(\kappa)$
and Injectivity Principle $IP(\kappa,\lambda)$ for regular $\kappa>\aleph_1$ and
$\lambda\leq\kappa$ which are formulated in terms of coloring the ordinals
$<\kappa$ by reals.}

\subsection{Long Borel Hierarchies}
We show that it is relatively consistent with ZF that the Borel hierarchy on
the reals has length $\omega_2$. This implies that $\omega_1$ has countable
cofinality, so the axiom of choice fails very badly in our model. A similar
argument produces models of ZF in which the Borel hierarchy has length any
given limit ordinal less than $\omega_2$, e.g., $\omega$ or
$\omega_1+\omega_1$.

\hfill{\texttt{http://www.math.wisc.edu/\~{}miller/res/index.html}}

\nby{Arnold W. Miller}

\nArxPaper{0705.0504}
{Rothberger's property in finite powers}
{Marion Scheepers}
{We show that several classical Ramseyan statements, and a forcing statement,
are each equivalent to having Rothberger's property in all finite powers.}

\nAMSPaper{http://www.ams.org/journal-getitem?pii=S0002-9939-07-08808-9}
{Special subsets of the reals and tree forcing notions}
{Marcin Kysiak, Andrzej Nowik, and Tomasz Weiss}
{We study relationships between classes of special subsets of the
reals (e.g., meager-additive sets, $\gamma$-sets, $C''$-sets, $\lambda$-sets) and the
ideals related to the forcing notions of Laver, Mathias, Miller
and Silver.}

\nArxPaper{0705.3085}
{All automorphisms of the Calkin algebra are inner}
{Ilijas Farah}
{We prove that it is relatively consistent with the usual axioms of
mathematics that all automorphisms of the Calkin algebra are inner. Together
with a 2006 Phillips--Weaver construction of an outer automorphism using the
Continuum Hypothesis, this gives a complete solution to a 1977 problem of
Brown--Douglas--Fillmore. We also give a simpler and self-contained proof of
the Phillips--Weaver result.}

\nArxPaper{0705.2867}
{Continuous selections and $\sigma$-spaces}
{Dusan Repovs, Boaz Tsaban, and Lyubomyr Zdomskyy}
{Assume that $X$ is a metrizable separable space, and each clopen-valued lower
semicontinuous multivalued map $\Phi$ from $X$ to $Q$ has a continuous selection. Our
main result is that in this case, $X$ is a $\sigma$-space. We also derive a partial
converse implication, and present a reformulation of the Scheepers Conjecture
in the language of continuous selections.}

\nArxPaper{0705.3877}
{On the closure of the diagonal of a $T_1$-space}
{Maria-Luisa Colasante and Dominic van der Zypen}
{Let $X$ be a topological space. The closure of $\Delta = \{(x, x) : x \in X\}$ in
$X \times X$ is a symmetric relation on $X$. We characterize those
equivalence relations on an infinite set that arise as the closure
of the diagonal with respect to a $T_1$-topology.}

\nArxPaper{0705.4297}
{Splitting families and the Noetherian type of $\beta\N\sm\N$}
{David Milovich}
{Extending some results of Malykhin, we prove several independence results
about base properties of $\beta\N\sm\N$ and its powers, especially the
Noetherian type $Nt(\beta\N\sm\N)$, the least $\kappa$ for which
$\beta\N\sm\N$ has a base that is $\kappa$-like with respect to
containment. For example, $Nt(\beta\N\sm\N)$ is never less than the
splitting number, but can consistently be that $\aleph_1$, $2^\alephes$,
$(2^\alephes)^+$, or strictly between $\aleph_1$ and $2^\alephes$.
$Nt(\beta\N\sm\N)$ is also consistently less than the additivity of the
meager ideal. $Nt(\beta\N\sm\N)$ is closely related to the existence of
special kinds of splitting families.
}

\nArxPaper{0706.0319}
{Even more simple cardinal invariants}
{Jakob Kellner}
{Using GCH, we force the following: There are continuum many simple cardinal
characteristics with pairwise different values.}

\nArxPaper{0706.1686}
{A classification of CO spaces which are continuous images of compact ordered spaces}
{Robert Bonnet and Matatyahu Rubin}
{A compact Hausdorff space $X$ is called a CO space, if every closed subset of $X$
is homeomorphic to an open subset of $X$. Every successor ordinal with its order
topology is a CO space. We find an explicit characterization of the class $K$ of
CO spaces which are a continuous image of a Dedekind complete totally ordered
set. (The topology of a totally ordered set is taken to be its order topology).
We show that every member of $K$ can be described as a finite disjoint sum of
very simple spaces. Every summand has either form: (1) $\mu + 1 + \nu^*$, where $\mu$
and $\nu$ are cardinals, and $\nu^*$ is the reverse order of $\nu$; or (2) the summand
is the 1-point-compactification of a discrete space with cardinality $\aleph_1$.}

\section{Problem of the Issue}

The general definitions of the properties considered here
are available in \cite{coc11}, but in our specific case
there are combinatorial reformulations \cite{BG}, which we
reproduce here for the reader's convenience.
For a $g\in\Bgp$, $f\in\NN$, and $n\in\N$,
$|g|\rest{[0,n)}\le f(n)$ means: For each $k\in[0,n)$ (that is, each $k<n$),
$|g(k)|\le f(n)$.

A subgroup $G$ of $\Bgp$ is \emph{Menger-bounded} if, and only if,
there is $f\in\NN$ such that
$$(\forall g\in G)(\exists^\oo n)\ |g|\rest{[0,n)}\le f(n).$$
$G^2$ is Menger-bounded if, and only if,
there is $f\in\NN$ such that
$$(\forall g_0,g_1\in G)(\exists^\oo n)(\forall i=0,1)\ |g_i|\rest{[0,n)}\le f(n).$$

In \cite{SqMen} we used a weak but unprovable hypothesis to prove
that there is a group $G\le\Bgp$ such that
$G$ is Menger-bounded, but $G^{2}$ is not Menger-bounded.

A subgroup $G$ of $\Bgp$ is \emph{Rothberger-bounded} if, and only if,
for each increasing $h\in\NN$, there is $\varphi:\N\to\FinSeqs{\Z}$
such that:
$$(\forall g\in G)(\exists n)\ g\rest{[0,h(n))} = \varphi(n).$$
$g\rest{[0,h(n))} = \varphi(n)$ means: For each $k\in[0,h(n))$ (that is, each $k<h(n)$),
$g(k)=\varphi(n)(k)$.

Clearly, each Rothberger-bounded group is Menger-bounded.

\begin{prob}
Does \textsf{CH} imply the existence of a group $G\le\Bgp$ such that $G$ is
Rothberger-bounded but $G^2$ is not Menger-bounded?
\end{prob}

This problem is raised in \cite{BG}, and was independently suggested by Ljub\v{s}a Ko\v{c}inac during
the SPM07 Workshop.

\ed